\documentclass{amsart}

\usepackage{amsfonts, amsmath, amssymb, amsthm, enumitem}
\usepackage[mathscr]{euscript}

\newcommand{\real}{\mathbb{R}}
\newcommand{\nat}{\mathbb{N}}
\newcommand{\dist}{\mathrm{dist}}

\newtheorem{thm}{Theorem}[section]
\newtheorem{cor}[thm]{Corollary}
\newtheorem{lem}[thm]{Lemma}

\theoremstyle{definition}

\theoremstyle{remark}
\newtheorem{rem}[thm]{Remark}
\newtheorem{ex}[thm]{Example}
\numberwithin{equation}{section}

\begin{document}

\title[Distortion from spheres into Euclidean spaces]{Distortion from spheres into Euclidean spaces}

\author[James Dibble]{James Dibble}
\address{Department of Mathematics and Statistics, University of Southern Maine, Portland, ME 04103, USA}
\email{james.dibble@maine.edu}


\begin{abstract}
Any function from a round $n$-dimensional sphere of radius $r$ into $n$-dimensional Euclidean space must distort the metric additively by at least $\displaystyle \frac{\pi r}{1 + \sqrt{1 - \frac{2}{n+2}}}$ if $n$ is even and $\displaystyle \frac{\pi r}{1 + \sqrt{1 - \frac{2(n+2)}{(n+1)(n+3)}}}$ if $n$ is odd. This is proved using a fixed-point theorem of Granas that generalizes the classical theorem of Borsuk--Ulam to set-valued functions.
\end{abstract}

\maketitle

\section{Introduction}

Let $X$ and $Y$ be metric spaces. The \textbf{distortion} of a nonempty relation $R$ from $X$ to $Y$ (i.e., a nonempty subset of $X \times Y$) is
\[
    \dist(R) = \sup \big\{ |d_Y(y_1,y_2) - d_X(x_1,x_2)| \,|\, (x_1,y_1),(x_2,y_2) \in R \big\},
\]
which takes values in $[0,\infty]$. In particular, if $f : X \to Y$ is a function, then
\[
    \dist(f) = \sup_{x_1,x_2 \in X} |d_Y(f(x_1),f(x_2)) - d_X(x_1,x_2)|.
\]
The idea of distortion appears in a number of places in geometry and analysis. For example, the Gromov--Hausdorff distance between two compact metric spaces is realized as half the infimum of the distortion of correspondences between them \cite[Theorem 7.3.25]{BuragoBuragoIvanov2001}.

Denote by $S^n_r$ the round $n$-dimensional sphere of radius $r$ in $\real^{n+1}$. That is,
\[
	S^n_r = \{ (x_1,\ldots,x_{n+1}) \in \real^{n+1} \,|\, x_1^2 + \cdots + x_{n+1}^2 = r^2 \}.
\]
The classical theorem of Borsuk--Ulam \cite{Borsuk1933} states that any continuous function from $S^n_r$ to $\real^n$ takes a pair of antipodal points to the same point and, consequently, has distortion at least $\pi r$. Since projection onto any hyperplane distorts by exactly that much, that is the minimum distortion over the set of continuous functions. The main results of this paper are a pair of more general lower bounds, depending on the parity of $n$, that hold over the set of all functions.

\begin{thm}\label{main theorem}
    Let $n \geq 1$. If $f : S^n_r \to \real^n$ is any function, then $\dist(f) \geq \displaystyle \frac{\pi r}{1 + \sqrt{1 - \frac{2}{n+2}}}$ if $n$ is even and $\dist(f) \geq \displaystyle \frac{\pi r}{1 + \sqrt{1 - \frac{2(n+2)}{(n+1)(n+3)}}}$ if $n$ is odd.
\end{thm}

\noindent Since both of the bounds in Theorem \ref{main theorem} are decreasing in $n$, they imply a universal lower bound for the distortion of functions from spheres into Euclidean spaces of no higher dimension.

\begin{cor}
	Let $n \geq 1$. If $f : S^n_r \to \real^m$ is any function and $m \leq n$, then $\dist(f) > \displaystyle \frac{\pi r}{2}$.
\end{cor}

\noindent The one-dimensional case of Theorem \ref{main theorem} was proved by the author and Elizabeth Kupin using a combinatorial argument. The general case is an application of a fixed-point theorem of Granas, which generalizes to upper semi-continuous set-valued functions the theorem of Borsuk--Ulam. Its proof requires the following upper bounds for the distance between the vertex sets of two simplices in $\real^n$ that intersect nontrivially.

\begin{thm}\label{vertex bound 2}
	Let $\Delta$ and $\Delta'$ be simplices in $\real^n$ that intersect nontrivially and have edges of length at most $L$. Then, there exist a vertex $v$ of $\Delta$ and a vertex $w$ of $\Delta'$ such that $d(v,w) \leq L \sqrt{1 - \frac{2}{n+2}}$ if $n$ is even and $d(v,w) \leq L\sqrt{1 - \frac{2(n+2)}{(n+1)(n+3)}}$ if $n$ is odd.
\end{thm}

\noindent Those bounds are sharp, as they are realized by two regular simplices of dimensions $\lfloor \frac{n}{2} \rfloor$ and $\lceil \frac{n}{2} \rceil$ that intersect orthogonally at their barycenters.\\

\noindent \textbf{Organization of the paper.} Section 2 contains an elementary combinatorial proof of the one-dimensional case of the main theorem. Section 3 contains background information about fixed-point theory for set-valued functions, including a theorem of Granas that generalizes the theorem of Borsuk--Ulam to set-valued functions. Section 4 contains information about the geometry of simplices in Euclidean spaces, leading to the proof of the upper bound for the distance between the vertex sets of two simplices that intersect nontrivially. Section 5 contains the proof of the main theorem.\\

\noindent \textbf{Acknowledgments.} The proof of the one-dimensional case was developed jointly with Elizabeth Kupin. The author is grateful to an anonymous referee for important corrections.

\section{The one-dimensional case}

Before proving the general case of Theorem \ref{main theorem}, it may help to give an elementary combinatorial proof for functions on $S^1_r$. This was developed jointly with Elizabeth Kupin. The quotient space $[0,2\pi] / \sim$, where $\sim$ identifies the endpoints, is homeomorphic to the circle $S^1_r$ by the map $C(\theta) = (r\cos \theta, r\sin \theta)$. Under this identification, consider the $m$ equally spaced points $x_k = C(2\pi k/m)$ on $S^1_r$, with their indices being taken modulo $m$. To simplify the exposition, $m$ will always be assumed to be odd.

\begin{lem}\label{one-dimensional lemma}
	Let $f : S^1_r \to \real$ be any function. Suppose that, for some $k$, one of the following holds:\\
	\indent (1) $f(x_k) \leq f(x_{k + \frac{m+1}{2}}) \leq f(x_{k+1})$\\
	\indent (2) $f(x_{k+1}) \leq f(x_{k+\frac{m+1}{2}}) \leq f(x_k)$\\
	Then, $\dist(f) \geq \frac{2\pi r(m-2)}{3m}$.
\end{lem}

\begin{proof}
	As neighbors on $S^1_r$, $x_k$ and $x_{k+1}$ are at a distance of $\frac{2\pi r}{m}$ from each other. At the same time, $x_k$ and $x_{k + \frac{m+1}{2}}$ are almost antipodal, i.e., the distance between them is $\frac{\pi r(m-1)}{m}$, which is also the distance between $x_{k+1}$ and $x_{k + \frac{m+1}{2}}$. Therefore,
	\begin{align*}
		\frac{2\pi r}{m} + \dist(f) &\geq d(f(x_k), f(x_{k+1}))\\
		&= d(f(x_k), f(x_{k + \frac{m+1}{2}})) + d(f(x_{k+\frac{m+1}{2}}), f(x_{k+1}))\\
		&\geq 2 \Big( \frac{\pi r(m-1)}{m} - \dist(f) \Big).
	\end{align*}
	The result follows by solving for $\dist(f)$.
\end{proof}

\begin{thm}
	If $f : S^1_r \to \real$ is any function, then $\dist(f) \geq 2\pi r/3$.
\end{thm}

\begin{proof}
	Let $m \geq 3$ be an odd integer, and let $x_k$ be the $m$ points in $S^1_r$ described above. Construct a graph on $\{ x_k \}$ by assigning an edge between $x_i$ and $x_j$ whenever $x_i$ and $x_j$ are almost antipodal, i.e., the distance between them is $\frac{\pi r(m-1)}{m}$. Since $m$ is odd, this graph is $2$-regular (i.e., each vertex is connected to exactly two others), and its edges form a cycle of length $m$. If, for some $m$, two almost antipodal points map to the same point under $f$, then the distortion of $f$ is at least $\frac{\pi r(m-1)}{m}$, and the result holds. If that never happens, then one may turn this into a directed graph by orienting each edge to point from $x_i$ to $x_j$ exactly when $f(x_i) < f(x_j)$. Since $m$ is odd, the edge cycle must contain a directed path of length at least two, which produces one of the two configurations in the hypothesis of Lemma \ref{one-dimensional lemma}. Therefore, $\dist(f) \geq \frac{2\pi r(m-2)}{3m}$. The result follows by letting $m \to \infty$.
\end{proof}

\noindent It's not difficult to show that the distortion bound is sharp in this case.

\begin{ex}
Under the identification of $S^1_r$ with $[0,2\pi] / \sim$, define a map $f : S^1_r \to [0,2\pi r/3)$ by $f(\theta) = r\theta/3$ for all $\theta \in [0,2\pi)$. Then, $\dist(f) = 2\pi r/3$.
\end{ex}

\section{Fixed-point theory for set-valued functions}

For any set $X$, denote by $\mathscr{P}(X)$ the power set of $X$. If $F : X \to \mathscr{P}(X)$ is a set-valued function, then $x \in X$ is a \textbf{fixed point} of $F$ whenever $x \in F(x)$. If $F : X \to \mathscr{P}(Y)$ and $A \subset X$, then the \textbf{image of $A$ under $F$} is $F(A) = \cup_{x \in A} F(x)$, and the \textbf{image of $F$} is $F(X)$. If $X$ and $Y$ are topological spaces, then $F : X \to \mathscr{P}(Y)$ is \textbf{upper semi-continuous} if, whenever $x_k$ and $y_k$ are sequences such that $x_k \to x$, $y_k \in F(x_k)$, and $y_k \to y$, it follows that $y \in F(x)$. When, in addition, $X$ is a metric space, a function $F : X \to \mathscr{P}(Y)$ is \textbf{compact} if the image of every bounded set in $X$ has compact closure in $Y$.

The first major result in this area was the celebrated fixed-point theorem of Kakutani \cite{Kakutani1941}, which generalizes Brouwer's fixed-point theorem \cite{Brouwer1911}.  If $E$ has an affine structure, denote by $CC(E)$ the set of nonempty, closed, and convex subsets of $E$.

\begin{thm}[Kakutani]
	Let $X \subset \real^n$ be nonempty, compact, and convex. If $F : X \to CC(X)$ is upper semi-continuous, then $F$ has a fixed point.
\end{thm}

\noindent This was generalized to nonempty, compact, and convex subsets of arbitrary Banach spaces by Bohnenblust--Karlin \cite{BohnenblustKarlin1950}.

In a similar vein, Granas \cite{Granas1959a} generalized the Leray--Schauder degree \cite{LeraySchauder1934} of a compact perturbation of the identity on a closed ball in a Banach space to the set-valued case. If $A$ is a subset of a Banach space $E$, a \textbf{compact perturbation of the identity on $A$} is a set-valued function $\phi : A \to CC(E)$ of the form $\phi(x) = x - \Phi(x)$ for a compact, upper semi-continuous function $\Phi : A \to CC(E)$. It is \textbf{nonvanishing} if $0 \not\in \phi(x)$ for all $x \in A$ or, equivalently, $\phi$ maps into $CC(E \setminus \{ 0 \})$. Note that $\phi$ is nonvanishing exactly when $\Phi$ has no fixed points. Granas defined the degree\footnote{Granas used the term ``characteristic'' for nonvanishing compact perturbations of the identity on spheres and reserved ``degree'' for compact perturbations of the identity on closed balls, with the degree of such a perturbation $g$ relative to a point $y_0$ not in the image of any point on the boundary sphere defined to be the characteristic of the restriction of $g - y_0$ to the boundary.} $\gamma(\phi)$ of a nonvanishing compact perturbation of the identity $\phi$ on a sphere $\mathscr{S} = \partial B(x,r) \subset E$ in such a way that it satisfies the following two properties  \cite[Th\'{e}or\`{e}me $\textrm{1}_\alpha$]{Granas1959a}:
\begin{enumerate}[label=\textbf{(\roman*)}, leftmargin=3em]
	\item If $\phi'$ is homotopic to $\phi$, in the sense that there exists $\psi : [0,1] \times \mathscr{S} \to CC(E \setminus \{ 0 \})$ of the form $\psi(t,x) = x - \Psi(t,x)$, where $\Psi$ is compact and upper semi-continuous, that agrees with $\phi'$ on $\{ 0 \} \times \mathscr{S}$ and with $\phi$ on $\{ 1 \} \times \mathscr{S}$, then $\gamma(\phi') = \gamma(\phi)$.
	\item If $f : S \to E$ is a continuous single-valued function such that $f(x) \in \phi(x)$ for all $x \in \mathscr{S}$, then $\gamma(\phi)$ agrees with the degree of $f$ defined in \cite{Rothe1938}. When $E = \real^n$, this is the usual degree of a map from an $(n-1)$-sphere into $\real^n \setminus \{ 0 \}$ due to Brouwer \cite{Brouwer1911}.
\end{enumerate}

\noindent With respect to this notion of degree, Granas in \cite[Theorem $\textrm{3}_\alpha$]{Granas1959b} proved the following generalization of a classical theorem of Borsuk \cite{Borsuk1933}.

\begin{thm}[Granas]\label{granas odd degree}
	Let $E$ be a Banach space, $S^\alpha$ the unit sphere in $E$, and $\phi$ a nonvanishing compact perturbation of the identity on $S^\alpha$. If $\phi$ is odd, i.e.,
\[
	\phi(-x) = -\phi(x) \textrm{ for all } x \in S^\alpha ,
\]
then $\gamma(\phi)$ is odd.
\end{thm}

\noindent In particular, if $\phi : S^n \to CC(\real^{n+1} \setminus \{ 0 \})$ is upper semi-continuous and odd, and if $\phi(S^n)$ has compact closure, then $\gamma(\phi)$ is odd.

The following corollary of Theorem \ref{granas odd degree} was not explicitly stated in \cite{Granas1959b}, but it almost surely was known at the time to be an immediate consequence. The proof is an adaptation to the set-valued case of a well-known argument that derives the theorem of Borsuk--Ulam from the single-valued version of Theorem \ref{granas odd degree}.

\begin{cor}\label{granas theorem}
	Let $F : S^n \to CC(\real^n)$ be an upper semi-continuous function such that $F(S^n)$ has compact closure. Then, there exists $x \in S^n$ such that $F(-x) \cap F(x) \neq \emptyset$.
\end{cor}

\begin{proof}
Assume the conclusion is false. Then, if one identifies $\real^n$ with the subset of $\real^{n+1}$ in which $x_{n+1} = 0$, the set-valued difference
\[
	\phi(x) = F(x) - F(-x) = \{ y - y' \,|\, y \in F(x) \textrm{ and } y' \in F(-x) \}
\]
is an upper semi-continuous odd function from $S^n$ into $CC(\real^{n+1} \setminus \{ 0 \})$ whose image has compact closure, so it has odd degree. However, for $p = (0,\ldots,0,-1)$, the compact, upper semi-continuous function
\[
	(t,x) \mapsto \{ (1-t)y + tp \,|\, y \in \phi(x) \}
\]
on $[0,1] \times S^n$ homotopically transforms $\phi$ into a constant single-valued map through upper semi-continuous compact set-valued functions $S^n \to CC(\real^{n+1} \setminus \{ 0 \})$. Thus, $\phi$ has degree zero, a contradiction.
\end{proof}

\noindent Corollary \ref{granas theorem} generalizes the theorem of Borsuk--Ulam to set-valued functions in the same way that Kakutani's theorem generalizes Brouwer's fixed-point theorem. 

\section{Simplicial geometry}

The goal of this section is to prove Theorem \ref{vertex bound 2}, which gives a sharp upper bound for the distance between the vertex sets of two simplices in $\real^n$ that intersect nontrivially and have edge lengths bounded above by a fixed constant. It seems reasonable to think that this result exists elsewhere in the literature, but a reference has proved elusive.

An \textbf{$m$-dimensional simplex} $\Delta(v_1,\ldots,v_{m+1})$ in $\real^n$ is the convex hull of $m+1$ affinely independent $v_1,\ldots,v_{m+1} \in \real^n$, i.e., points such that the set $\{ v_1 - v_{m+1}, \ldots, v_m - v_{m+1} \}$ is linearly independent. As a set,
\[
	\Delta(v_1,\ldots,v_{m+1}) = \Big\{ \sum_{i=1}^{m+1} t_i v_i \,\big|\, t_1,\ldots,t_{m+1} \geq 0 \textrm{ and } \sum_{i=1}^{m+1} t_i = 1 \Big\}.
\]
A \textbf{$k$-dimensional face} of $\Delta(v_1,\ldots,v_{m+1})$ is a $k$-dimensional simplex obtained as the convex hull of exactly $k + 1$ of the $v_i$. A $k$-dimensional face of an $m$-dimensional simplex is \textbf{proper} if $k < m$. The \textbf{vertices} of $\Delta(v_1,\ldots,v_{m+1})$ are its $0$-dimensional faces, and its \textbf{edges} are its $1$-dimensional faces. In other words, its vertices are the $v_i$, and its edges are the line segments connecting them. A simplex is \textbf{regular} if its edges all have the same length, in which case all of its faces are also regular. The \textbf{barycenter} of $\Delta(v_1,\ldots,v_{m+1})$ is $\sum_{i=1}^{m+1} v_i/(m+1)$. The \textbf{boundary} $\partial \Delta$ of a simplex $\Delta$ is the union of its proper faces, and its \textbf{interior} $\Delta^\circ$ is $\Delta \setminus \partial \Delta$. These agree with its relative boundary and relative interior, respectively, as a topological subspace of $\real^n$.

The \textbf{Chebyshev radius} of a compact subset $X$ of $\real^n$ is the minimum $r \geq 0$ such that $X$ is contained in a closed ball of radius $r$, and the \textbf{Chebyshev center} of $X$ is the center of the smallest such ball. By the strict convexity of the distance function, every compact set in $\real^n$ has a unique Chebyshev center: if it had Chebyshev radius $R$ and two different Chebyshev centers, then a ball of radius less than $R$ around their midpoint would cover the set.

The following classical result of Jung \cite{Jung1901} will be needed.

\begin{thm}[Jung]\label{jung's theorem}
	If $X$ is a compact subset of $\real^n$ with diameter $D$, then $X$ is contained in a closed ball of radius at most $D\sqrt{\frac{n}{2(n+1)}}$.
\end{thm}

\noindent It follows that the Chebyshev radius of a $k$-dimensional simplex with edges of length at most $L$ is bounded above by $L\sqrt{\frac{k}{2(k+1)}}$. This is realized by regular simplices, whose Chebyshev centers are their barycenters.

\begin{lem}\label{vertex distance bound}
	Let $\Delta_k = \Delta(v_1,\ldots,v_{k+1})$ and $\Delta_m = \Delta(w_1,\ldots,w_{m+1})$ be simplices in $\real^n$ that intersect nontrivially and have Chebyshev radiuses $R_v$ and $R_w$, respectively. Then, there exist $i$ and $j$ such that $d(v_i,w_j) \leq \sqrt{R_v^2 + R_w^2}$.
\end{lem}

\begin{proof}
	Let $x = \sum_{i=1}^{k+1} \alpha_i v_i = \sum_{j=1}^{m+1} \beta_j w_j$ be a point of intersection of $\Delta_k$ and $\Delta_m$, where $0 \leq \alpha_i,\beta_j \leq 1$ and $\sum_{i=1}^{k+1} \alpha_i = \sum_{j=1}^{m+1} \beta_j = 1$. Without loss of generality, after applying a translation, one may suppose that $x$ is the origin, so that
\[
	0 = |x|^2 = x \cdot x = \big( \sum_{i=1}^{k+1} \alpha_i v_i \big) \cdot \big( \sum_{j=1}^{m+1} \beta_j w_j \big) = \sum_{i=1}^{k+1} \sum_{j=1}^{m+1} \alpha_i \beta_j (v_i \cdot w_j).
\]
Denote by $p_v$ and $p_w$ the respective Chebyshev centers of $\Delta_k$ and $\Delta_m$. Then,
\[
	R_v^2 \geq |v_i - p_v|^2 = |v_i|^2 + |p_v|^2 - 2(v_i \cdot p_v)
\]
for each $i$. Thus,
\begin{equation}\label{Chebyshev radius bound 1}
\begin{aligned}
	R_v^2 &= \sum_{i=1}^{k+1} \sum_{j=1}^{m+1} \alpha_i \beta_j R_v^2 \geq \sum_{i=1}^{k+1} \sum_{j=1}^{m+1} \alpha_i \beta_j [|v_i|^2 + |p_v|^2 - 2(v_i \cdot p_v)]\\
	&= \sum_{i=1}^{k+1} \alpha_i |v_i|^2 + |p_v|^2 - 2 \big[ \big( \underbrace{\sum_{i=1}^{k+1} \alpha_i v_i}_{x} \big) \cdot p_v \big]\\
	&= \sum_{i=1}^{k+1} \alpha_i |v_i|^2 + |p_v|^2 .
\end{aligned}
\end{equation}
Similarly,
\begin{equation}\label{Chebyshev radius bound 2}
	R_w^2 \geq \sum_{j=1}^{m+1} \beta_j |w_j|^2 + |p_w|^2 .
\end{equation}
Assume that $|v_i - w_j| > \sqrt{R_v^2 + R_w^2}$ for all $i$ and $j$. Then,
\[
	|v_i|^2 + |w_j|^2 - 2(v_i \cdot w_j) = |v_i - w_j|^2 > R_v^2 + R_w^2 ,
\]
so
\begin{align*}
	R_v^2 + R_w^2 &= \sum_{i=1}^{k+1} \sum_{j=1}^{m+1} \alpha_i \beta_j (R_v^2 + R_w^2) < \sum_{i=1}^{k+1} \sum_{j=1}^{m+1} \alpha_i \beta_j [|v_i|^2 + |w_j|^2 - 2(v_i \cdot w_j)]\\
	&= \sum_{i=1}^{k+1} \alpha_i |v_i|^2 + \sum_{j=1}^{m+1} \beta_j |w_j|^2 - 2 \Big[ \sum_{i=1}^{k+1} \sum_{j=1}^{m+1} \alpha_i \beta_j (v_i \cdot w_j) \Big]\\
	&= \sum_{i=1}^{k+1} \alpha_i |v_i|^2 + \sum_{j=1}^{m+1} \beta_j |w_j|^2 .
\end{align*}
Combining this with \eqref{Chebyshev radius bound 1} and \eqref{Chebyshev radius bound 2} shows that $|p_v|^2 + |p_w|^2 < 0$, a contradiction.
\end{proof}

\begin{rem}
	The case $m = 0$ of Lemma \ref{vertex distance bound} states that every point of $\Delta_k$ is at a distance of no more than $R_v$ from at least one vertex of $\Delta_k$.
\end{rem}

\noindent If regular $k$- and $m$-dimensional simplices $\Delta_k$ and $\Delta_m$ with edges of length $L$ intersect orthogonally at their barycenters, then the distance between the vertices of $\Delta_k$ and those of $\Delta_m$ is
\[
L \sqrt{\frac{2km + k + m}{2(k+1)(m+1)}} = L \sqrt{1 - \frac{1}{2}\Big( \frac{1}{k+1} + \frac{1}{m+1} \Big)}.
\]
It follows from Jung's theorem and Lemma \ref{vertex distance bound} that this quantity is an upper bound for the distance between the vertex sets of more general simplices that intersect nontrivially.

\begin{lem}\label{vertex bound 1}
	Let $\Delta_k = \Delta(v_1,\ldots,v_{k+1})$ and $\Delta_m = \Delta(w_1,\ldots,w_{m+1})$ be $k$- and $m$-dimensional simplices, respectively, with edges of length at most $L$. If $\Delta_k$ and $\Delta_m$ intersect nontrivially, then there exist $i$ and $j$ such that
\[
	d(v_i,w_j) \leq L\sqrt{1 - \frac{1}{2}\Big( \frac{1}{k+1} + \frac{1}{m+1} \Big)}.
\]
\end{lem}

\noindent It remains in this section to show that, if $\Delta_k$ and $\Delta_m$ are simplices in $\real^n$ that intersect nontrivially and have edge lengths bounded above by $L$, then the largest possible distance between their vertex sets is given by the bound in Lemma \ref{vertex bound 1} when $k = \lfloor \frac{n}{2} \rfloor$ and $m = \lceil \frac{n}{2} \rceil$.

\begin{lem}\label{intersect proper face}
	Let $\Delta_k = \Delta(v_1,\ldots,v_{k+1})$ and $\Delta_m = \Delta(w_1,\ldots,w_{m+1})$ be simplices in $\real^n$ that intersect nontrivially. If $k + m > n$, then a proper face of $\Delta_k$ intersects $\Delta_m$ or vice versa.
\end{lem}

\begin{proof}
	Since $\Delta_k$ and $\Delta_m$ intersect nontrivially, so do their respective affine hulls $H_k$ and $H_m$. Since $k + m > n$, $H_k \cap H_m$ has positive dimension. Let $x \in \Delta_k \cap \Delta_m$, and let $u \in \real^n$ be any unit vector such that $x + u \in H_k \cap H_m$. Since the line $t \mapsto x + t u$ must eventually leave the compact sets $\Delta_k$ and $\Delta_m$, there are points at which it hits $\partial \Delta_k$ and points at which it hits $\partial \Delta_m$. Any nearest such point to $x$ satisfies the conclusion of the lemma.
\end{proof}

\begin{proof}[Proof of Theorem \ref{vertex bound 2}]
Suppose that $\Delta_k$ and $\Delta_m$ are simplices in $\real^n$ that intersect nontrivially and have edges of length at most $L$. If $k + m > n$, then, by Lemma \ref{intersect proper face}, the hypotheses of the theorem will still hold if either $\Delta_k$ or $\Delta_m$ is replaced with one of its proper faces. Therefore, without loss of generality, one may suppose that $k + m \leq n$. In this case,
\[
	\sqrt{1 - \frac{1}{2}\Big( \frac{1}{k+1} + \frac{1}{m+1} \Big)} \leq \sqrt{1 - \frac{n+2}{2(k+1)(n-k+1)}}.
\]
Note that
\[
	(k+1)(n-k+1) = - \Big( k-\frac{n}{2} \Big)^2 + \frac{n^2}{4} + n + 1.
\]
When $n$ is even, this is maximized at $k = \frac{n}{2}$; when $n$ is odd, it is maximized at $k = \frac{n-1}{2}$ and $k = \frac{n+1}{2}$. The result follows from Lemma \ref{vertex bound 1} by substituting those values of $k$ into the right-hand side of the above inequality.
\end{proof}

\section{Proof of the main theorem}

Let $f : X \to E$ be any function from a metric space $X$ into a normed vector space $E$. Define a set-valued function $F : X \to CC(E)$ by
\[
	F(x) = \bigcap_{\varepsilon > 0} \overline{\mathrm{conv}\big( f(B(x,\varepsilon)) \big)},
\]
where $\mathrm{conv}(Y)$ and $\overline{Y}$ denote the convex hull and closure, respectively, of a set $Y$ and $B(x,\varepsilon)$ denotes the open ball of radius $\varepsilon$ around $x$. This function captures the ``convex essence'' of $f$ at each point. Since $f(x) \in f(B(x,\varepsilon))$ for all $\varepsilon > 0$, $F(x) \neq \emptyset$. As the intersection of closed and convex sets, each $F(x)$ is, in fact, closed and convex.

\begin{lem}\label{upper semi-continuous}
	For any $f : X \to E$, the function $F$ defined above is upper semi-continuous.
\end{lem}

\begin{proof}
	Suppose $x_k \in X$, $y_k \in F(x_k)$, $x_k \to x$, and $y_k \to y$. Let $\varepsilon > 0$, and fix $K \in \nat$ such that $x_k \in B(x,\varepsilon)$ for all $k \geq K$. For each such $k$, set $r_k = \varepsilon - d(x_k,x)$. From the definition of $F(x_k)$, one has that
\[
	y_k \in \overline{\mathrm{conv} \big( f(B(x_k, r_k)) \big)} .
\]
By the triangle inequality, $B(x_k,r_k) \subseteq B(x,\varepsilon)$, so $y_k \in \overline{\mathrm{conv} \big( f(B(x,\varepsilon)) \big)}$, and, consequently, $y \in \overline{\mathrm{conv} \big( f(B(x,\varepsilon)) \big)}$. Since the choice of $\varepsilon > 0$ was arbitrary, $y \in F(x)$.
\end{proof}

\noindent It will be helpful to prove an alternate characterization of each $F(x)$, but this requires a condition on $f$. A function $f : X \to \real^n$ is \textbf{locally bounded} if each point of $X$ has a neighborhood whose image under $f$ is bounded. In what follows, a theorem of Carath\'{e}odory \cite{Caratheodory1911} (proved in full generality by Steinitz \cite{Steinitz1913}) will be used repeatedly. Carath\'{e}odory's theorem states that, for $Y \subseteq \real^n$, $\mathrm{conv}(Y)$ is the union of all simplices of dimension at most $n$ with vertices drawn from $Y$. This may be used to show that $\overline{\mathrm{conv}(S)} = \mathrm{conv}(\overline{S})$ for bounded $S$. In particular, $\mathrm{conv}(Y)$ is compact whenever $Y$ is compact.

\begin{rem}
	While it is true that $\mathrm{conv}(\overline{S}) \subseteq \overline{\mathrm{conv}(S)}$ for all $S$, the reverse inclusion might not hold if $S$ is unbounded, as shown by $S = \{ (0,0) \} \cup \{ (1/n,n) \,|\, n \in \nat \}$.
\end{rem}

\begin{lem}\label{nested bounded subsets}
	If $\{ Y_t \,|\, 0 < t < T \}$ is a nested family of bounded subsets of $\real^n$ (i.e., $Y_s \subseteq Y_t$ whenever $s \leq t$), then $\mathrm{conv}(\bigcap_{0 < t < T} \overline{Y_t}) = \bigcap_{0 < t < T} \overline{ \mathrm{conv}(Y_t) }$.
\end{lem}

\begin{proof}
	It will first be shown that $\mathrm{conv}(\bigcap_{0 < t < T} Z_t) = \bigcap_{0 < t < T} \mathrm{conv}(Z_t)$ for a nested family of compact $Z_t$. To see the inclusion $\mathrm{conv} ( \bigcap_{0 < t < T} Z_t ) \subseteq \bigcap_{0 < t < T} \mathrm{conv}(Z_t)$, it suffices by Carath\'{e}odory's theorem to show that, for all $v_1,\ldots,v_\ell \in \bigcap_{0 < t < T} Z_t$, $\mathrm{conv} \{ v_1,\ldots,v_\ell \} \subseteq \bigcap_{0 < t < T} \mathrm{conv}(Z_t)$. This latter inclusion follows from the fact that, for each $0 < t < T$, $v_i \in Z_t$ for each $i$, so $\mathrm{conv} \{ v_1,\ldots,v_\ell \} \subseteq \mathrm{conv}(Z_t)$.

	To see the inclusion $\bigcap_{0 < t < T} \mathrm{conv}(Z_t) \subseteq \mathrm{conv} ( \bigcap_{0 < t < T} Z_t )$, choose an arbitrary $z \in \bigcap_{0 < t < T} \mathrm{conv}(Z_t)$. Since $z \in \mathrm{conv} (Z_{1/m})$ for each sufficiently large $m$, there exist, by Carath\'{e}odory's theorem, $v_m^1,\ldots,v_m^{n+1} \in Z_{1/m}$ with the property that $z \in \mathrm{conv} \{ v_m^1,\ldots,v_m^{n+1} \}$. Since the $Z_{1/m}$ are compact and nested, one may, by passing to a subsequence, suppose without loss of generality that $v_m^i \to v_i \in Z_{1/M}$ for all sufficiently large $M$. Since the $Z_t$ are nested, $v_i \in \bigcap_{0 < t < T} Z_t$. By writing $z$ as convex combinations of the $v_m^i$ and, if necessary, passing to another subsequence so that the coefficients converge, one finds that $z \in \mathrm{conv} \{ v_1,\ldots,v_{n+1} \} \subseteq \mathrm{conv} ( \bigcap_{0 < t < T} Z_t )$.

	In the general case, since the $Y_t$ are bounded, each $\overline{Y_t}$ is compact. Thus,
\[
	\mathrm{conv} \big( \bigcap_{0 < t < T} \overline{Y_t} \big) = \bigcap_{0 < t < T} \mathrm{conv}(\overline{Y_t}) = \bigcap_{0 < t < T} \overline{ \mathrm{conv}(Y_t) } .
\]
\end{proof}

\noindent For each $x \in X$, the sets $f(B(x,\varepsilon))$ are nested. Thus, the following is an immediate consequence of Lemma \ref{nested bounded subsets}.

\begin{lem}\label{rewrite F}
	If $f : X \to \real^n$ is locally bounded, then
\[
	F(x) = \mathrm{conv} \big( \bigcap_{\varepsilon > 0} \overline{f(B(x,\varepsilon))} \big).
\]
\end{lem}

\begin{proof}[Proof of Theorem \ref{main theorem}]
	Let $f : S^n_r \to \real^n$ be any function. If $f(S^n_r)$ does not have compact closure, then $f$ is unbounded. In this case, $\dist(f) = \infty$, and the result holds. If $f(S^n_r)$ has compact closure, then, by Lemma \ref{rewrite F},
\[
	F(x) =  \mathrm{conv} \big( \bigcap_{\varepsilon > 0} \overline{f(B(x,\varepsilon))} \big).
\]
By Lemma \ref{upper semi-continuous}, $F$ is upper semi-continuous, so, by Corollary \ref{granas theorem}, there exists $x \in S_r^n$ such that $F(x) \cap F(-x) \neq \emptyset$. By Carath\'{e}odory's theorem, there exist affinely independent $v_1,\ldots,v_{k+1} \in \bigcap_{\varepsilon > 0} \overline{f(B(x,\varepsilon))}$ and $w_1,\ldots,w_{m+1} \in \bigcap_{\varepsilon > 0} \overline{f(B(-x,\varepsilon))}$ such that the simplices $\Delta(v_1,\ldots,v_{k+1})$ and $\Delta(w_1,\ldots,w_{m+1})$ intersect nontrivially.

Let
\[
	L = \max \{ d(v_i,v_{i'}), d(w_j,w_{j'}) \,|\, 1 \leq i,i' \leq k+1 \textrm{ and } 1 \leq j,j' \leq m+1 \}.
\]
For each $\varepsilon > 0$, $v_i \in \overline{f(B(x,\varepsilon))}$ and $w_j \in \overline{f(B(-x,\varepsilon))}$ for all $i$ and $j$. Suppose $d(v_i,v_{i'}) = L$, and choose $y_\ell^i,y_\ell^{i'} \in B(x,\varepsilon)$ such that $f(y_\ell^i) \to v_i$ and $f(y_\ell^{i'}) \to v_{i'}$ as $\ell \to \infty$. By the triangle inequality,
\[
	\dist(f) + 2\varepsilon > \dist(f) + d(y_\ell^i,y_\ell^{i'}) \geq d(f(y_\ell^i),f(y_\ell^{i'})) \to d(v_i,v_{i'}) = L,
\]
so $L \leq \dist(f) + 2\varepsilon$. If $d(w_j,w_{j'}) = L$, one obtains the same inequality. By a similar argument, for each $i$ and $j$, $d(v_i,w_j) \geq \pi r - \dist(f) - 2\varepsilon$. Since these hold for each $\varepsilon > 0$, $L \leq \dist(f)$ and $d(v_i,w_j) \geq \pi r - \dist(f)$. By Theorem \ref{vertex bound 2}, there exist $i_0$ and $j_0$ such that $d(v_{i_0},w_{j_0}) \leq L\sqrt{1 - \frac{2}{n+2}}$ if $n$ is even and $d(v_{i_0},w_{j_0}) \leq L\sqrt{1 - \frac{2(n+2)}{(n+1)(n+3)}}$ if $n$ is odd. Thus,
\[
	\pi r - \dist(f) \leq d(v_{i_0},w_{j_0}) \leq \dist(f) \sqrt{1 - \frac{2}{n+2}}
\]
if $n$ is even, and
\[
	\pi r - \dist(f) \leq d(v_{i_0},w_{j_0}) \leq \dist(f) \sqrt{1 - \frac{2(n+2)}{(n+1)(n+3)}}
\]
if $n$ is odd. The result follows by solving each of the above for $\dist(f)$.
\end{proof}

\bibliography{bibliography}
\bibliographystyle{amsplain}

\end{document}